\documentclass[12pt,a4paper,draft]{article}
\usepackage[catalan,spanish,english]{babel}
\usepackage[psamsfonts]{amssymb}
\usepackage{cmmib57}
\usepackage{exscale}
\usepackage{amsmath}
\usepackage{amscd}
\usepackage{latexsym}
\usepackage{graphicx}
\usepackage{pdfsync}
\newtheorem{thm}{Theorem}[section]

\newtheorem{prop}[thm]{Proposition}

\newtheorem{rem}{Remark}[section]
\newcommand{\qed}{\quad{$\square$}} 

\newcommand {\IR}{\mathbb{R}}

\newcommand {\F}{\mbox{${\mathcal F}$}}

\newcommand {\hac}{\mathcal H}
\newcommand{\hact}{\mathcal{H}_T}

\newcommand{\ep}{\varepsilon}

\newcommand{\tbar}{\bar{t}}

\newcommand{\beq}{\begin{equation}}
\newcommand{\eeq}{\end{equation}}
\newcommand{\beqn}{\begin{equation*}}
\newcommand{\eeqn}{\end{equation*}}

\newcommand{\tf}{\mathcal{F}}

\begin{document}

\begin{titlepage}
\null \vspace{0.5cm}
\begin{center}
{\Large\bf  A Laplace principle\\[1mm]
 for a stochastic wave equation\\[1mm]
in spatial dimension three\\[2mm]}
\medskip

by\\
\vspace{7mm}

\begin{tabular}{l@{\hspace{10mm}}l@{\hspace{10mm}}l}
{\sc V\'{\i}ctor Ortiz-L\'opez}$\,^{(\ast)}$ &and &{\sc Marta Sanz-Sol\'e}$\,^{(\ast)}$\\
{\small vortize@ub.edu}         &&{\small marta.sanz@ub.edu}\\
&&{\small http://www.mat.ub.es/$\sim$sanz}\\
\end{tabular}
\begin{center}
{\small Facultat de Matem\`atiques}\\
{\small Universitat de Barcelona } \\
{\small Gran Via de les Corts Catalanes, 585} \\
{\small E-08007 Barcelona, Spain} \\
\end{center}
\end{center}

\vspace{1cm}

\noindent{\bf Abstract:} We consider a stochastic wave equation in spatial dimension three, driven by 
a Gaussian noise, white in time and with a stationary spatial covariance. The free terms
are nonlinear with Lipschitz continuous coefficients. Under suitable conditions on the 
covariance measure, Dalang and Sanz-Sol\'e \cite{dss09} have proved the existence of a random
field solution with H\"older continuous sample paths, jointly in both arguments, time and 
space. By perturbing the driving noise with a multiplicative parameter $\varepsilon\in]0,1]$,
 a family of probability laws corresponding to the respective solutions to the equation is obtained.
Using the weak convergence approach to large deviations developed in \cite{de97}, we prove that
this family  satisfies a Laplace principle in the H\"older norm.

\medskip


\vspace{1 cm}

\noindent


{\begin{itemize} \item[$^{(\ast)}$] Supported by the grant MTM
2009-07203 from the \textit{Direcci\'on General de
Investigaci\'on, Ministerio de Ciencia e Innovaci\'on, Spain.}
\end{itemize}}
\end{titlepage}

\newpage

 \section{Introduction}
 \label{s0}

 We consider the stochastic wave equation in spatial dimension $d=3$
 \begin{equation} 
 \label{0.1}
\begin{cases}\left(\frac{\partial^2}{\partial t^2} - \varDelta \right)u(t,x) = \sigma\big(u(t,x)\big)\dot F(t,x) 
+ b \big(u(t,x)\big),\ t\in ]0,T],\\ 
u(0,x) = v_0(x),\\
 \frac{\partial}{\partial t}u(0,x) = \tilde v_0(x),  \ x\in\IR^3,
 \end{cases}
 \end{equation}
where 
$\varDelta$ denotes the  Laplacian on $\IR^3$. The coefficients $\sigma$ and $b$ are Lipschitz continuous functions and
the  process $\dot F$ is the {\it formal} derivative of a Gaussian random field, 
white in time and correlated in space. More precisely, for any $d\ge 1$,
let $\mathcal{D}(\mathbb{R}^{d+1})$ be the space of Schwartz test functions 
 and let $\varGamma$ be a non-negative and 
non-negative definite tempered measure on $\IR^d$. 
Then, on some probability space, there exists a Gaussian process  $F = \left(F(\varphi),\ \varphi\in \mathcal{D}(\mathbb{R}^{d+1})\right)$
 with mean zero and covariance functional  
\begin{equation}\label{cov1}
   E\big(F(\varphi) F(\psi)\big) = \int_{\mathbb{R}_+} ds \int_{\IR^d} \varGamma(dx) (\varphi(s)*\tilde\psi(s))(x), 
\end{equation}
where $\tilde\psi(s)(x)=\psi(s)(-x)$ and the notation ``$\ast$" means the convolution operator. As has been proved in \cite{dalangfrangos}, the process $F$ can be extended to
a martingale measure $M=\left(M_t(A), t\ge 0, A\in\mathcal{B}_b(\IR^d)\right)$, where $\mathcal{B}_b(\IR^d)$ denotes the set of bounded Borel sets of $\IR^d$.

For any $\varphi, \psi\in \mathcal{D}(\IR^n)$, define the inner product
\beqn
\langle \varphi,\psi\rangle_{\mathcal{H}}=\int_{\IR^n}\varGamma(dx)(\varphi\ast\tilde\psi)(x)
\eeqn
and denote by $\mathcal{H}$ the Hilbert space obtained by the completion of $\mathcal{D}(\IR^n)$ with the inner product 
$\langle \cdot,\cdot\rangle_{\mathcal {H}}$. Using the theory of stochastic integration with respect to martingale measures
(see for instance \cite{walsh}), the stochastic integral $B_t(h):=\int_0^tds \int_{\IR^d} h(y) M(ds,dy)$ is well defined, and for any $h\in\hac$
with $\Vert h\Vert_\hac=1$, the process $(B_t(h), t\in[0,T])$ is a standard Wiener process.
In addition, for any fixed $t\in[0,T]$, the mapping $h\to B_t(h)$ is linear. Thus, the process $(B_t, t\in[0,T])$ is a cylindrical Wiener process on
$\mathcal{H}$ (see \cite{dz} for a definition of this notion). Let $(e_k, k\ge 1)$ be a complete orthonormal system of $\mathcal{H}$. Clearly, $B_k(t):=\int_0^t ds \int_{\IR^d} e_k(y) M(ds,dy)$, $k\ge1$,
defines a sequence of independent, standard Wiener processes and we have the representation 
\beq
\label{0.7}
B_t=\sum_{k\ge 1} B_k(t) e_k.
\eeq
Let $\mathcal{F}_t$, $t\in[0,T]$, be the $\sigma$-field generated by the random variables $(B_k(s), s\in[0,t], k\ge 1)$.  $(\mathcal{F}_t)$-predictable processes
$\Phi\in L^2(\Omega\times [0,T]; \hac)$ can be integrated with respect to the cylindrical Wiener process $(B_t, t\in[0,T])$ and the stochastic integral 
$\int_0^t \Phi(s) dB_t$ coincides with the It\^o stochastic integral with respect to the infinite dimensional Brownian motion $(B_k(t), t\in [0,T], k\ge 1)$,
$ \sum_{k\ge 1}\int_0^t \langle\Phi(s),e_k\rangle_\hac dB_k(t)$.

We shall consider the mild formulation of equation (\ref{0.1}),
\begin{align} \label{0.2}
u(t,x)&= w(t,x)+\sum_{k\ge 1}\int_0^t \langle G(t-s,x-\cdot) \sigma(u(s,\cdot)),e_k\rangle_\hac dB_k(s)\nonumber\\
& +\int_0^t  [G(t-s)\ast b(u(s,\cdot))](x) ds,
\end{align}
$t\in[0,T]$, $x\in \IR^3$. Here 
\beqn
w(t,x)= \left(\frac{d}{dt}G(t)\ast v_0\right)(x)+(G(t)\ast\tilde v_0)(x), 
\eeqn
and $G(t)=\frac{1}{4\pi t}\sigma_t$, where $\sigma_t$ denotes the uniform surface measure (with total mass
$4\pi t^2$) on the sphere of radius $t$.  

\smallskip

Throughout the paper, we will consider the following set of assumptions.
\smallskip

\noindent{\bf(H)}
\begin{enumerate}
\item The coefficients $\sigma$, $b$ are real Lipschitz continuous functions.
\item The spatial covariance measure $\varGamma$ is absolutely continuous with respect to Lebesgue measure and the density is
$f(x)= \varphi(x)|x|^{-\beta}$, $x\in\IR^3\backslash\{0\}$. The function $\varphi$ is bounded and positive, $\varphi\in\mathcal{C}^1(\IR^3)$, 
$\nabla \varphi\in\mathcal{C}_b^\delta(\IR^3)$
(the space of bounded and H\"older continuous functions with exponent $\delta\in]0,1]$) and $\beta\in]0,2[$. 
\item The initial values $v_0$, $\tilde v_0$ are bounded and such that $v_0\in\mathcal{C}^2(\IR^3)$, $\nabla v_0$ is bounded and $\varDelta v_0$ and $\tilde v_0$ are H\"older continuous with degrees $\gamma_1,
\gamma_2\in]0,1]$, respectively.
\end{enumerate}
We remark that the assumptions on $\varGamma$ imply
\beq
\label{noise}
\sup_{t\in[0,T]}\int_{\IR^3} \vert \F(G(t))(\xi)\vert^2 \mu(d\xi)<\infty,
\eeq
where $\tf$ denotes the Fourier transform operator and $\mu=\tf^{-1}\varGamma$. This is a relevant condition in connection with the definition of the stochastic integral with respect to the martingale measure $M$ (\cite{dalang}).

The set of hypotheses {\bf (H)} are used in Chapter 4 of \cite{dss09} to prove a theorem on existence and uniqueness of solution to equation (\ref{0.2}) and the properties of the sample paths. More precisely, under a slightly weaker set of assumptions than {\bf (H)} (not requiring boundedness of the functions $v_0$, $\tilde v_0$, $\nabla v_0$),  Theorem 4.11 in \cite{dss09} states that
 for any $q\in[2,\infty[$, $\alpha\in]0,\gamma_1\wedge\gamma_2\wedge \frac{2-\beta}{2}\wedge\frac{1+\delta}{2}[$, there exists $C>0$ such that for 
$(t,x),(\bar t,y)\in[0,T]\times D$,
\beq
\label{0.3}
E(|u(t,x)-u(\bar t,y)|^q) \le C(|t-\tbar|+|x-y|)^{\alpha q},
\eeq
where $D$ is a fixed bounded domain of $\IR^3$.
Consequently, a.s., the stochastic process $(u(t,x), (t,x)\in[0,T]\times D)$ solution of (\ref{0.2}) has $\alpha$-H\"older continuous sample paths, jointly in $(t,x)$.

The reason for strengthening the assumptions of \cite{dss09} is to ensure that
\beq
\label{initial}
\sup_{(t,x)\in[0,T]\times \IR^3}\vert w(t,x)\vert  < \infty
\eeq
(see Hypothesis 4.1 and Lemma 4.2 in \cite{dqs09}), a condition that is needed in the proof of Theorem \ref{t2} below.
This is in addition to (4.19) in \cite{dss09}, which provides an estimate of a fractional Sobolev norm of the function $w$.

We notice that in \cite{dss09}, the mild formulation of equation (\ref{0.1}) is stated using the stochastic integral developed in \cite{dm03}. Recent results by Dalang and Quer-Sardanyons
(see \cite{dqs09}, Proposition 2.11 and Proposition 2.6 (b)) show that this formulation is equivalent to (\ref{0.2}).

In this paper, we consider the family 
of stochastic wave equations
\begin{align} \label{0.4}
u^\ep(t,x)&= w(t,x)+\sqrt \ep\sum_{k\ge 1}\int_0^t \langle G(t-s,x-\cdot) \sigma(u(s,\cdot)),e_k\rangle_\hac dB_k(s)\nonumber\\
& +\int_0^t  [G(t-s)\ast b(u^\ep(s,\cdot))](x) ds,
\end{align}
$\ep\in]0,1]$, and we establish a large deviation principle for the family $(u^\ep, \ep\in]0,1])$ in a Polish space closely related to $\mathcal{C}^{\alpha}([0,T]\times D)$,
the space of functions defined on $[0,T]\times D$, H\"older continuous jointly in its two arguments, of degree $\alpha\in \mathcal{I}$, where 
\beqn
\mathcal{I}:=\left]0,\gamma_1\wedge\gamma_2\wedge \frac{2-\beta}{2}\wedge\frac{1+\delta}{2}\right[.
\eeqn

To formulate the large deviation principle, we should consider a Polish
space carrying the probability laws of the family $(u^\ep, \ep>0)$. This cannot
be $\mathcal{C}^{\alpha}([0,T]\times D)$, since this space is not separable. 
Instead, we consider the space $\mathcal{C}^{\alpha^\prime,0}([0,T]\times D)$ of
H\"older continuous functions $g$ of degree $\alpha^\prime<\alpha$, with 
modulus of continuity
\beqn
O_g(\delta):=\sup_{|t-s|+|x-y|<\delta}\frac{|g(t,s)-g(s,y)|}{(|t-s|+|x-y|)^{\alpha^\prime}}
\eeqn
satisfying $\lim_{\delta\to 0^+}O_g(\delta)=0$. This is a Banach space and 
$\mathcal{C}^{\alpha}([0,T]\times D)\subset \mathcal{C}^{\alpha^\prime,0}([0,T]\times D)$.

In the sequel, we shall denote by $(\mathcal{E}_\alpha, \Vert\cdot\Vert_\alpha)$ the Banach space $\mathcal{C}^{\alpha,0}([0,T]\times D)$
endowed with the H\"older norm of degree $\alpha$,
and consider values of $\alpha\in\mathcal{I}$.
\medskip

Let $\hact=L^2([0,T]; \hac)$. For any $h\in\hact$, we consider the deterministic evolution equation 
\begin{align}
\label{0.5}
V^h(t,x)&= w(t,x)\nonumber\\
&+\int_0^t  \langle G(t-s,x-\cdot) \sigma(V^h(s,\cdot)),h(s,\cdot)\rangle_{\hac} ds \nonumber\\
&+\int_0^t \left[G(t-s)\ast b(V^h(s)\right](x)\ ds.
\end{align}
The second term on the right-hand side of this equation can be written as
\beqn
\sum_{k\ge 1}\int_0^t \langle G(t-s,x-\cdot) \sigma(V^h(s,\cdot)),e_k\rangle_{\hac} h_k(s)\ ds,
\eeqn
with $h_k(t)=\langle h(t),e_k\rangle_\hac$, $t\in[0,T]$, $k\ge 1$.

Existence and uniqueness of solution of equation (\ref{0.5}) can be proved in a similar (but easier) way than for (\ref{0.2}). 
This will be obtained in the next section as a by-product of Theorem \ref{t2}, where it is also proved that
 $V^h\in \mathcal{E}_\alpha$.
We will denote by $\mathcal{G}^0: \hact\longrightarrow \mathcal{E}_\alpha$ the mapping defined by  $\mathcal{G}^0(h)=V^h$.

For any $f\in \mathcal{E}_\alpha$ define
\beq
\label{0.6}
I(f)=\inf_{h\in\hact: \mathcal{G}^0(h)=f}\left\{\frac{1}{2}\Vert h\Vert_{\hact}^2\right\}
\eeq
and for any $A\subset \mathcal{E}_\alpha$, $I(A)=\inf\{I(f), f\in A\}$.
\medskip

The main result of this paper is the following theorem.
\begin{thm}
\label{t1}
Assume that the set of hypotheses {\bf(H)} are satisfied. Then, the family $\{u^\ep, \ep\in]0,1]\}$ given by (\ref{0.4}) satisfies a large deviation principle on $\mathcal{E}_\alpha$ with rate function $I$ given by \eqref{0.6}. That means,
for any closed subset $F\in\mathcal{E}_\alpha$ and any open subset $G\in \mathcal{E}_\alpha$,
\begin{align*}
&\limsup_{\ep \to 0^+}\ep \log P(u^\ep\in F)\le -I(F),\\
&\liminf _{\ep \to 0^+}\ep \log P(u^\ep\in G)\ge -I(G).
\end{align*}
\end{thm}

In the proof of this theorem, we will use the weak convergence approach to large deviations developed in \cite{de97}. 
An essential ingredient of this method is a variational representation for a reference Gaussian process (Brownian motion when studying diffusion processes, or different generalizations
of infinite-dimensional Wiener process when dealing with stochastic partial differential equations).
As it is shown in \cite{bd00}, a variational representation for an infinite-dimensional Brownian motion along with a transfer principle based on compactness and weak convergence, allow to derive a large deviation principle for some functionals of this process. This method has been applied in \cite{bdm08} to establish a large deviation principle to reaction-diffusion systems considered in \cite{k92} and also in several subsequent papers, for instance in \cite{ss06}, \cite{dm09}, \cite{z09}.
We next give the ingredients
for the proof of Theorem \ref{t1} based on this method. 
\medskip

\noindent {\sl Variational representation of infinite dimensional Brownian motion} 
\smallskip

Let $B=(B_k(t), t\in[0,T], k\ge 1)$ be a sequence of independent standard Brownian motions.
Denote by $\mathcal{P}(l^2)$ the set of predictable processes belonging to $L^2(\Omega\times [0,T]); l^2)$ and let $g$ be a real-valued, bounded,
Borel measurable function defined on $\mathcal{C}([0,T]; \IR^\infty)$. Then,
\beq
\label{0.8}
-\log E(\exp[-g(B)]) = \inf_{u\in \mathcal{P}(l^2)} E\left(\frac{1}{2} \Vert u\Vert_{L^2([0,T]; l^2)}^2 + g\left(B+\int_0^\cdot u\right)\right)
\eeq
( see Theorem 2 in \cite{bdm08}).
\medskip

\noindent{\sl Weak regularity}

Denote by $\mathcal{P}_\hac$ the set of predictable processes belonging to $L^2(\Omega\times [0,T]); \hac)$. For any $N>0$, we define 
\begin{align*}
\hact^N&=\{h\in\hact: \Vert h\Vert_{\hact} \le N\},\\
\mathcal{P}_{\hac}^N&=\{v\in \mathcal{P}_{\hac}: v\in\hact^N, a.s.\},
\end{align*}
and we consider $\hact^N$ endowed with the weak topology of $\hact$.

For any $v\in \mathcal{P}_\hac^N$, $\ep\in]0,1]$, let $u^{\ep,v}$ be the solution to
\begin{align} \label{0.9}
u^{\ep,v}(t,x)&= w(t,x)+\sqrt \ep\sum_{k\ge 1}\int_0^t\langle G(t-s,x-\cdot) \sigma(u^{\ep,v}(s,\cdot)),e_k\rangle_\hac d B_k(s)\nonumber\\
& + \int_0^t \langle G(t-s,x-\cdot) \sigma(u^{\ep,v}(s,\cdot)),v(s,\cdot)\rangle_\hac\ ds \nonumber\\
& +\int_0^t  [G(t-s)\ast b(u^{\ep,v}(s,\cdot))](x) ds.
\end{align}
We will prove in Theorem \ref{t2} that this equation has a unique solution and that $u^{\ep,v}\in\mathcal{E}_\alpha$ with
 $\alpha\in \mathcal{I}$.
 
Consider the following conditions:
\begin{description}
\item {(a)} The set $\{V^h, h\in\hact^N\}$ is a compact subset of $\mathcal{E}_\alpha$, where $V^h$ is the solution of (\ref{0.5}).
\item {(b)} For any family $(v^\ep, \ep>0)\subset \mathcal{P}_\hac^N$ which converges in distribution as $\ep\to 0$ to $v\in\mathcal{P}_\hac^N$, as
$\hact^N$-valued random variables, we have
\beqn
\lim_{\ep\to 0}u^{\ep,v^\ep} = V^v,
\eeqn
in distribution, as $\mathcal{E}_\alpha$-valued random variables. 
\end{description}
Here $V^v$ stands for the solution of (\ref{0.5}) corresponding to a $\hact^N$-valued random variable $v$ (instead of a deterministic function $h$). The solution is a stochastic process $\{V^h(t,x), (t,x)\in[0,T]\times \IR^3\}$ defined path-wise by (\ref{0.5}).

According to \cite{bdm08}, Theorem 6 applied to the functional  $\mathcal{G}: \mathcal{C}([0,T]; \IR^\infty)\to \mathcal{E}_\alpha$, $\mathcal{G}(\sqrt \ep B):=u^\ep$ (the solution of (\ref{0.4})),  and
$\mathcal{G}^0: \hact\to \mathcal{E}_\alpha$, $\mathcal{G}^0(h):= V^h$ (the solution of (\ref{0.5})), conditions (a) and (b) above imply the validity of Theorem \ref{t1}.

\section{Laplace principle for the wave equation}
\label{s1}

Following the discussion of the preceding section, the proof of Theorem \ref{t1} will consist of checking that conditions (a) and (b) above hold true. As we next show, both conditions will follow from a single continuity result. Indeed, 
the set $\hact^N$ is a compact subset of $\hact$ endowed with the weak topology (see \cite{lax}, Chapter 12, Theorem 4). Thus, (a)
can be obtained by proving that the mapping $h\in\hact^N\mapsto V^h\in \mathcal{E}_\alpha$ is continuous with respect to the weak
topology. For this, we consider a sequence $(h_n, n\ge 1)\subset \hact^N$ and $h\in\hact^N$ satisfying $\lim_{n\to\infty}\Vert h_n-h\Vert_w=0$, which means that for any $g\in\hact$, 
$\lim_{n\to\infty}\langle h_n-h,g\rangle_{\hact}=0$, and 
we will prove that
\beq
\label{1.1}
\lim_{n\to\infty}\Vert V^{h_n}-V^h\Vert_{\alpha}=0.
\eeq
As for (b), we invoke Skorohod Representation Theorem and rephrase this condition as follows. On some probability space $(\bar\Omega,\bar\tf, \bar P)$,
consider a sequence of independent Brownian motions $\bar B=\{\bar B_k, k\ge 1\}$ along with the corresponding filtration $(\bar\tf_t, t\in[0,T])$, where $\bar\tf_t$ is the $\sigma$-field generated by the random variables $(\bar B_k(s), s\in[0,t],k\ge1)$. Furthermore, consider
a family of $(\bar\tf_t)$-predictable processes $(\bar v^\ep, \ep>0, \bar v)$ belonging to $L^2(\bar\Omega\times [0,T]; \hac)$ taking values on $\hact^N$, $\bar P$ a.s.,
such that the joint law of $(v^\ep, v, B)$ (under $P$) coincides with that of $(\bar v^\ep, \bar v, \bar B)$ (under $\bar P$) and such that, 
\beqn
\lim_{\ep\to 0}\Vert \bar v^\ep - \bar v \Vert_w=0, \ \bar P- a.s.
\eeqn
as $\hact^N$-valued random variables.
Let $\bar u^{\ep,\bar v^\ep}$ be the solution to a similar equation as (\ref{0.9}) obtained by changing $v$ into $\bar v^\ep$ and $B_k$ into $\bar B_k$. Then, we will prove that for any $q\in[0,\infty[$,
\beq
\label{1}
\lim_{\ep\to 0}\bar E\left(\left\Vert \bar u^{\ep,\bar v^\ep}-V^{\bar v}\right\Vert_\alpha^q\right)=0,
\eeq
where $\bar E$ denotes the expectation operator on $(\bar\Omega,\bar\tf, \bar P)$.
Notice that, if in (\ref{0.9}) we consider $\ep=0$ and $v:=h\in \mathcal{P}_{\hac}^N$ deterministic, we obtain the equation satisfied by $V^h$. Consequently, the convergence (\ref{1.1}) can be obtained as a particular case of (\ref{1}).

Therefore, we will focus our efforts on the proof of (\ref{1}).
In the sequel, we shall omit any reference to the bars in the notation, for the sake of simplicity. 

Accoding to Lemma A1 in \cite{bmss95}, the proof of (\ref{1}) can be carried out into two steps:
\begin{enumerate}
\item {\sl Estimates on increments}
\begin{align}
&\sup_{\ep\ge 1}E\left(\left\vert \left[u^{\ep,v^\ep}(t,x) - V^v(t,x)\right]-\left[u^{\ep,v^\ep}(r,z)-V^v(r,z)\right]\right\vert^q\right)\nonumber\\
&\quad \le C [\vert t-r\vert + \vert x-z\vert]^{\alpha q}.\label{2}
\end{align}
\item {\sl Pointwise convergence}
\beq
\label{p}
\lim_{\ep\to 0}E\left(\vert u^{\ep,v^\ep}(t,x) - V^v(t,x)\vert^q\right)=0.
\eeq
\end{enumerate}
Here, $q\in[1,\infty[$, $(t,x), (r,z)\in [0,T]\times D$ and $\alpha\in \mathcal{I}$.

Before proving these facts, we will address the problem of giving a rigorous formulation of (\ref{0.9}). As we have already mentioned, the stochastic integral with respect to $(B_k, k\ge 1)$ in (\ref{0.9}) is equivalent to the stochastic integral $\int_0^t\int_{\IR^3} G(t-s,x-y)\sigma(u^{\ep,v^\ep}(s,y)) M(ds,dy)$ considered in the sense of \cite{dm03}. We recall that such an integral is defined for stochastic processes $Z=(Z(s,\cdot), s\in[0,T])$ with values in $L^2(\IR^3)$ a.s., adapted and mean-square continuous,
 and the integral
 \beq
 \label{idm}
v_{G,Z}^t(\star):= \sum_{k\ge 1}\int_0^t  \langle G(t-s,\star-\cdot) Z(s,\cdot),e_k(s,\cdot)\rangle_\hac dB_k(s)
\eeq
satisfies
\beq
\label{3}
E\left(\Vert v_{G,Z}^t\Vert^2_{L^2(\IR^3)}\right)=
\int_0^tds \int_{\IR^3} d\xi E(\vert \mathcal{F}Z(s)(\xi)\vert^2)\int_{\IR^3} \mu(d\eta)\vert\mathcal{F}G(t-s)(\xi-\eta)\vert^2.
\eeq
(see \cite{dm03}, Theorem 6).

As a function of the argument $x$, for any $v\in\mathcal{P}_{\hac}^N$, the path-wise integral 
\beqn
\int_0^t \langle G(t-s,x-\cdot) \sigma(u^{\ep,v}(s,\cdot)),v(s,\cdot)\rangle_\hac\ ds,
\eeqn
is also a well-defined $L^2(\IR^3)$-valued random variable. 
Indeed, let $Z$
be a stochastic process satisfying the hypotheses described before. Set 
\beq
\label{pidm}
\nu_{G,Z}^t(\star):=\int_0^t  \langle G(t-s,\star-\cdot) Z(s,\cdot),v(s,\cdot)\rangle_\hac\ ds.
\eeq
By Cauchy-Schwarz' inequality applied to the inner product on $\hact$, we have
\begin{align*}
 \Vert \nu_{G,Z}^t\Vert^2_{L^2(\IR^3)}&\le N^2
\int_{\IR^3} dx \int_0^t ds \Vert G(t-s,x-\cdot) Z(s,\cdot)\Vert^2_\hac\\
& = N^2 \int_0^tds \int_{\IR^3} d\xi \vert \mathcal{F}Z(s)(\xi)\vert^2\int_{\IR^3} \mu(d\eta)\vert\mathcal{F}G(t-s)(\xi-\eta)\vert^2,
\end{align*}
where the last equality is derived following the arguments for the proof of Theorem 6 in \cite{dm03}. We recall that this formula is firstly established for $Z$ sufficiently smooth and by smoothing $G$ by convolution with an approximation of the identity. The extension of the formula to the standing assumptions is done by a limit procedure.

From this, we clearly have
\begin{align}
\label{4}
& E\left(\Vert \nu_{G,Z}^t\Vert^2_{L^2(\IR^3)}\right)\le N^2\nonumber\\
& \quad\times\int_0^tds \int_{\IR^3} d\xi E(\vert \mathcal{F}Z(s)(\xi)\vert^2)\int_{\IR^3} \mu(d\eta)\vert\mathcal{F}G(t-s)(\xi-\eta)\vert^2.
\end{align}
\begin{rem}
\label{r1}
Up to a positive constant, the $L^2(\Omega;  L^2(\IR^3))$-norm of the stochastic integral $v^t_{G,Z}$ and the path-wise integral $\nu^t_{G,Z}$ are bounded by the same expression.
\end{rem}
Let $\mathcal{O}$ be a bounded or unbounded open subset of $\IR^3$, $q\in [1,\infty[$, $\gamma\in]0,1[$. We denote by $W^{\gamma,q}(\mathcal{O})$ the fractional Sobolev Banach space consisting of functions $\varphi: \IR^3 \rightarrow \IR$ such that
\beqn
\Vert \varphi\Vert_{W^{\gamma,q}(\mathcal{O})}:=\left(\Vert \varphi\Vert_{L^q(\mathcal{O})}^q
+ \Vert \varphi\Vert^q_{\gamma,q,\mathcal{O}}\right)^{\frac{1}{q}}<\infty,
\eeqn
where
\beqn
\Vert \varphi\Vert_{\gamma,q,\mathcal{O}}=\left(\int_{\mathcal{O}} dx \int_{\mathcal{O}} dy\frac{\vert \varphi(x)-\varphi(y)\vert^q}{\vert x-y\vert^{3+\gamma q}}\right)^{\frac{1}{q}}.
\eeqn 
For any $\ep>0$, we denote by $\mathcal{O}^\ep$ the $\ep$-enlargement of of $\mathcal{O}$, that is,
\beqn
\mathcal{O}^\ep=\{x\in \IR^3: d(x,\mathcal{O})<\ep\}.
\eeqn

In the proof of Theorem \ref{t2} below, we will use a smoothed version of the fundamental solution $G$ defined as follows.  
Consider a function $\psi\in\mathcal{C}^\infty(\IR^3; \IR_+)$ with support included in the unit ball, such that $\int_{\IR^3} \psi(x) dx=1$. For any $t\in]0,1]$ and $n\ge 1$, set
\beqn
\psi_n(t,x)=\left(\frac{n}{t}\right)^3\psi\left(\frac{n}{t}x\right),
\eeqn
and
\beq
\label{1.7}
G_n(t,x)=\left(\psi_n(t,\cdot) \ast G(t)\right)(x).
\eeq
Notice that, for any $t\in[0,T]$, ${\text{supp}}\ G_n(t,\cdot)\subset B_{t(1+\frac{1}{n})}(0)$.
\begin{rem}
\label{r2}
Since $G_n(t)$ is smooth and has compact support, $v_{G_n,z}^t(x)$ is well-defined
as a Walsh  stochastic integral and this integral defines a random field indexed by $(t,x)$. By Burkholder's inequality, for any $q\in[2,\infty[$,
\beqn
E\left(\vert v_{G_n,z}^t(x)\vert^q\right) \le C E\left(\int_0^t  \Vert G(t-s, x-\cdot) Z(s,\cdot)\Vert_{\hac}^2\ ds\right)^{\frac{q}{2}}.
\eeqn
As for the path-wise integral $\nu_{G_n,z}^t(x)$, by applying Cauchy-Schwarz' inequality to the inner product on $\hact$, we have
\beqn
E\left(\vert \nu_{G_n,z}^t(x)\vert^q\right) \le N^q E\left(\int_0^t \Vert G(t-s, x-\cdot) Z(s,\cdot)\Vert_{\hac}^2\ ds\right)^{\frac{q}{2}}.
\eeqn
Hence, as in Remark \ref{r1}, up to a constant, $L^q(\Omega)$-estimates for both type of integrals at fixed $(t,x)\in[0,T]\times \IR^3$ coincide. 
\end{rem}

The following proposition is the analogue of Theorem 3.1 in \cite{dss09} for the path-wise integral $\nu_{G,Z}^t$.

\begin{prop}
\label{p4}
Fix $q\in]3,\infty[$ and a bounded domain $\mathcal{O}\subset \IR^3$. Supose that
\beqn
\tau_q(\beta,\delta):=\left(\frac{2-\beta}{2}\wedge\frac{1+\delta}{2}\right)-\frac{3}{q}>0
\eeqn
and fix $\gamma\in]0,1[$, $\rho\in]0,\tau_q(\beta,\delta)\wedge \gamma[$.
Let $\{Z_t, t\in[0,T]\}$ be a $L^2(\IR^3)$-valued, $(\tf_t)$-adapted, mean-square continuous stochastic process. Assume that for some fixed $t\in[0,T]$,
\beqn
\int_0^t  E\left(\Vert Z(s)\Vert^q_{W^{\gamma,q}(\mathcal{O}^{t-s})}\right)\ ds < \infty.
\eeqn
We have the following estimates:
\begin{align}
&E\left(\Vert \nu_{G,Z}^t\Vert^q_{L^q(\mathcal{O})}\right) \le C \int_0^t  E\left(\Vert Z(s)\Vert^q_{L^q(\mathcal{O}^{t-s})}\right)\ ds,\label{1.8}\\
&E\left(\Vert \nu_{G_n,Z}^t\Vert^q_{\rho,q,\mathcal{O}}\right) \le C \int_0^t  E\left(\Vert Z(s)\Vert^q_{W^{\rho,q}(\mathcal{O}^{(t-s)(1+\frac{1}{n})})}\right)\ ds,\label{1.9}\\
&E\left(\Vert \nu_{G,Z}^t\Vert^q_{\rho,q,\mathcal{O}}\right) \le C \int_0^t  E\left(\Vert Z(s)\Vert^q_{W^{\rho,q}(\mathcal{O}^{(t-s)})}\right)\ ds\label{1.10}.
\end{align}
Consequently,
\beq
\label{1.11}
E\left(\Vert \nu_{G,Z}^t\Vert^q_{W^{\rho,q}(\mathcal{O})}\right) \le C \int_0^t  E\left(\Vert Z(s)\Vert^q_{W^{\rho,q}(\mathcal{O}^{(t-s)})}\right)\ ds.
\eeq
\end{prop}
\noindent{\it Proof}. By virtue of Remark \ref{r2}, we see that the estimate \eqref{1.8} follows from the same arguments used in \cite{dss09}, Proposition 3.4. We recall that this proposition is devoted to prove an analogue property for the stochastic integral $v_{G,Z}^t$. In the very same way, \eqref{1.9} is established using the arguments of the proof of Proposition 3.5 in \cite{dss09}.
Then, as in \cite{dss09},  \eqref{1.10} is obtained from \eqref{1.9} by applying Fatou's lemma. Finally, \eqref{1.11} is a consequence of \eqref{1.8}, \eqref{1.10} and the definition of the fractional Sobolev norm $\Vert \cdot\Vert_{W^{\rho,q}(\mathcal{O})}$.

\hfill\qed
\bigskip

Next, we present an analogue of Theorem 3.8 \cite{dss09} for the path-wise integral $\nu^t_{G,Z}$, which gives the sample path properties in the argument $t$ for this integral. As in Proposition \ref{p4}, $\mathcal{O}$ is a bounded domain in $\IR^3$.
\begin{prop}
\label{p5}
Consider a stochastic process $\{Z_t, t\in[0,T]\}$, $(\tf_t)$-adapted, with values in $L^2(\IR^3)$, mean-square continuous. Assume that for some fixed $q\in]3,\infty[$ and
$\gamma\in]\frac{3}{q},1[$,
\beqn
\sup_{t\in[0,T]} E\left(\Vert Z(t)\Vert^q_{W^{\gamma,q}(\mathcal{O}^{T-t})}\right) <\infty.
\eeqn
Then the stochastic process $\{\nu^t_{G,Z}(x), t\in[0,T]\}$, $x\in\mathcal{O}$, satisfies
\beq
\label{time}
\sup_{x\in \mathcal{O}}E\left(\vert \nu^t_{G,Z}(x)-\nu^{\bar t}_{G,Z}(x)\vert^{\bar q}\right) \le C \vert t-\bar t\vert^{\rho\bar q},
\eeq
for each $t, \bar t\in[0,T]$, any $\bar q\in]2,q[$, $\rho\in ]0,(\gamma-\frac{3}{q})\wedge(\frac{2-\beta}{2})\wedge(\frac{1+\delta}{2})[$.
\end{prop}
\noindent{\it Proof}. We follow the same scheme as in the proof of \cite{dss09}, Theorem 3.8. To start with, we should prove an analogue of (3.26) in \cite{dss09}, with $v^{\bar t}_{G_n,Z}$,  $v^{t}_{G_n,Z}$ replaced by $\nu^{\bar t}_{G_n,Z}$,  $\nu^{t}_{G_n,Z}$, respectively. Once again, we apply Remark \ref{r2}, obtaining similar upper bounds for the $L^{\bar q}(\Omega)$-moments (up to a positive constant) as for the stochastic integrals considered in the above mentioned reference. More precisely, assume $0\le t<\bar t\le T$; by applying Cauchy-Schwarz' inequality to the inner product on $\hact$, we obtain 
\begin{align*}
&E\left(\left\vert\int_t^{\bar t}  \langle G_n(\bar t-s,x-\cdot) Z(s,\cdot), v(s,\cdot)\rangle_\hac\ ds\right\vert^{\bar q}\right)\\
&\quad \le N^{\bar q} E\left(\int_0^{\bar t-t}  \left\Vert  G_n(s, x-\cdot)Z(\bar t-s, \cdot)\right\Vert^2_{\hac}\ ds\right)^{\frac{\bar q}{2}},\\
&E\left(\left\vert\int_0^{t}  \langle (G_n(\bar t-s,x-\cdot)-G_n(t-s,x-\cdot) Z(s,\cdot), v(s,\cdot)\rangle_\hac\ ds\right\vert^{\bar q}\right)\\
&\quad \le N^{\bar q} E\left(\int_0^{t}  \left\Vert (G_n(\bar t-s,x-\cdot)-G_n(t-s,x-\cdot) Z(s,\cdot)\right\Vert^2_\hac\ ds\right)^{\frac{\bar q}{2}}.
\end{align*}
These are, up to a positive constant, the same upper bounds obtained in \cite{dss09} for the expressions termed $T_1^n(t,\bar t,x)$ and $T_2^n(t,\bar t,x)$, respectively.
After this remark, the proof follows the same arguments as in \cite{dss09}.

\hfill\qed

For any $t\in[0,T]$, $a\ge 1$, let $K^D_a(t)=\{y\in\IR^3: d(y,D)\le a(T-t)\}$. For $a=1$, we shall simply write $K^D(t)$; this is the light cone of $\{T\}\times D$.  
\medskip

In the next theorem, the statement on existence and uniqueness of solution, as well as (\ref{t2.1}), extend Theorem 4.3 in \cite{dqs09}, while (\ref{t2.2}) and
\eqref{t2.3} are extensions of the inequality (4.24) of Theorem 4.6  and (4.41) of Theorem 4.11 in \cite{dss09}, respectively. Indeed in the cited references, the results apply to Equation \eqref{0.2} while in the next theorem, they apply to Equation \eqref{0.9}.

\begin{thm}
\label{t2}
Assuming {\bf (H)}, the following statements hold true:

There exists a unique random field solution to (\ref{0.9}), $\{u^{\ep,v}(t,x), (t,x)\in[0,T]\times \IR^3\}$, and this solution satisfies
\begin{align}
&\sup_{\ep\in]0,1], v\in\mathcal{P}_\hac^N}\sup_{(t,x)\in[0,T]\times \IR^3}E\left(\vert u^{\ep,v}(t,x)\vert^q\right) < \infty, \label{t2.1}\\
&\sup_{\ep\in]0,1], v\in\mathcal{P}_\hac^N}\sup_{t\in[0,T]}E\left(\Vert u^{\ep,v}(t)\Vert^q_{W^{\alpha,q}(K^D(t))}\right) < \infty, \label{t2.2}
\end{align}
for any $q\in[2,\infty[$, $\alpha\in\mathcal{I}$.

Moreover, for any $q\in[2,\infty[$ and $\alpha\in\mathcal{I}$, there exists $C>0$ such that for 
$(t,x),(\bar t,y)\in[0,T]\times D$,
\beq
\label{t2.3}
\sup_{\ep\in]0,1], v\in\mathcal{P}_\hac^N}E(|u^{\ep,v}(t,x)-u^{\ep,v}(\bar t,y)|^q) \le C(|t-\bar t|+|x-y|)^{\alpha q}.
\eeq

Thus, a.s., $\{u^{\ep,v}(t,x), (t,x)\in [0,T]\times D\}$ has H\"older continuous sample paths
of degree $\alpha\in\mathcal{I}$, jointly in $(t,x)$.
\end{thm}

\noindent{\it Proof}. 
For the sake of simplicity, we shall consider $\ep=1$ and write $u^v$ instead of $u^{\ep,v}$.

We start by proving existence and uniqueness along with (\ref{t2.1}).  For this, we will follow the method of the proof of \cite{dqs09}, Theorem 4.3 (borrowed from \cite{mss99}, Theorem 1.2 and \cite{dalang}, Theorem 13). It is based on the Picard iteration scheme:
\begin{align}
u^{v,(0)}(t,x)& = w(t,x),\nonumber\\
u^{v,(n+1)}(t,x)&= w(t,x) + \sum_{k\ge 1}\int_0^t\langle G(t-s,x-\cdot) \sigma(u^{v,(n)}(s,\cdot)),e_k\rangle_\hac d B_k(s)\nonumber\\
& + \int_0^t  \langle G(t-s,x-\cdot) \sigma(u^{v,(n)}(s,\cdot)),v(s,\cdot)\rangle_\hac\ ds\nonumber\\
& +\int_0^t  [G(t-s)\ast b(u^{v,(n)}(s,\cdot))](x)\ ds,\quad n\ge 0.\label{picard}
\end{align}
The steps of the proof are as follows. Firstly, we check that
\beq
\label{s1bis2}
\sup_{v\in\mathcal{P}_\hac^N}\sup_{(t,x)\in[0,T]\times \IR^3}E\left(\vert u^{v,(n)}(t,x)\vert^q\right)<\infty,
\eeq
and then 
\beq
\label{s1bis}
\sup_{n\ge 0}\sup_{v\in\mathcal{P}_\hac^N}\sup_{(t,x)\in[0,T]\times \IR^3}E\left(\vert u^{v,(n)}(t,x)\vert^q\right)<\infty.
\eeq
Secondly, by setting
\beqn
M_n(t):=\sup_{(s,x)\in[0,t]\times \IR^3} E\left(\vert u^{v,(n+1)}(s,x)-u^{v,(n)}(s,x)\vert^q\right),\ n\ge 0,
\eeqn
we prove
\beq
\label{s2}
M_n(t)\le C \int_0^t  M_{n-1}(s)\left(1+\int_{\IR^3} \mu(d\xi)\vert \tf G(t-s)(\xi)\vert^2\right)\ ds.
\eeq
With these facts, we conclude that $(u^{v,(n)}(t,x), n\ge 0)$ converges uniformly in $(t,x)$ in $L^q(\Omega)$ to a limit $u^{v}(t,x)$ which satisfies equation (\ref{0.9}) with $\ep=1$.

In comparison with the proof of Theorem 4.3 in \cite{dqs09}, establishing \eqref{s1bis2}--\eqref{s2} requires additionally the analysis of the term given by the path-wise integral
\beq
\label{s3}
\mathcal{I}^{v,(n+1)}:=\int_0^t \langle G(t-s,x-\cdot) \sigma(u^{v,(n)}(s,\cdot)),v(s,\cdot)\rangle_\hac\ ds, \ n\ge 0.
\eeq
This is done as follows. We assume that \eqref{s1bis2} holds true for some $n\ge 0$. This is definitely the case for $n=0$ (see \eqref{initial}).
By applying Cauchy-Schwarz inequality on the Hilbert space $\hact$, and since $\Vert v\Vert_{\hact}\le N$ a.s., we have
\begin{align*}
&E\left(\left\vert \int_0^t \langle G(t-s,x-\cdot) \sigma(u^{v,(n)}(s,\cdot)),v(s,\cdot)\rangle_\hac\ ds \right\vert^q\right)\\
&\le N^q E\left(\int_0^t  \left\Vert G(t-s,x-\cdot) \sigma(u^{v,(n)}(s,\cdot))\right\Vert_\hac^2 \ ds \right)^{\frac{q}{2}}.
\end{align*}

Notice that, by applying Burkholder's inequality to the stochastic integral term in \eqref{picard}, we obtain
\begin{align*}
&E\left(\left\vert \sum_{k\ge 1}\int_0^t\langle G(t-s,x-\cdot) \sigma(u^{v,(n)}(s,\cdot)),e_k\rangle_\hac d B_k(s)\right\vert^q\right)\\
&\le C E\left( \int_0^t \left\Vert G(t-s,x-\cdot) \sigma(u^{v,(n)}(s,\cdot))\right\Vert_\hac^2\ ds \right)^{\frac{q}{2}}.
\end{align*}
Thus, as has already been mentioned in Remark \ref{r2}, up to a positive constant, $L^q(\Omega)$ estimates of the stochastic integral and of the path-wise integral $\mathcal{I}^{v,(n+1)}$ lead to the same  upper bounds. 

This simple but important remark yields  the extension of properties \eqref{s1bis2}--\eqref{s2}, which are valid for Equation \eqref{0.2}, as is proved in Theorem 4.3 in \cite{dqs09}, to Equation \eqref{0.9}
 with $\ep=1$ and actually, for any $\ep\in]0,1]$. In fact, those properties can be proved to hold uniformly in $\ep\in]0,1]$.
 \medskip
 
 Let us now argue on the validity of \eqref{t2.2}. We will follow the programme of Section 4.2 in \cite{dss09}, taking into account the new term 
 \beqn
 \int_0^t \langle G(t-s,x-\cdot) \sigma(u^{v}(s,\cdot)),v(s,\cdot)\rangle_\hac\ ds 
 \eeqn
 of Equation \eqref{0.9} (with $\ep=1$) that did not appear in \cite{dss09}. This consists of the following steps. 
 
 Firstly, we need an extension of Proposition 4.3 in \cite{dss09}. This refers to an approximation of the localized version of \eqref{0.9} on a light cone. In the approximating sequence, the fundamental solution $G$ of the wave equation is replaced by a the smoothed version $G_n$ defined in \eqref{1.7}.
 Going through the proof of that Proposition, we see that for the required extension the term
 \beqn
   M_n(t):= E\left(\Vert v^t_{G_n,Z}-v^t_{G,Z}\Vert^q_{L^q(K_a^D(t))}\right),
   \eeqn
   with $Z(s,y)= \sigma(u^v(s,y))1_{K_a^D(s)}(y)$, 
 should be replaced by  
  \beqn
   \tilde M_n(t):= E\left(\Vert v^t_{G_n,Z}-v^t_{G,Z}\Vert^q_{L^q(K_a^D(t))}\right)
   +E\left(\Vert \nu^t_{G_n,Z}-\nu^t_{G,Z}\Vert^q_{L^q(K_a^D(t))}\right),
   \eeqn
   where we have used the notation introduced in \eqref{idm}, \eqref{pidm}.
 Then we should prove that $\lim_{n\to\infty} \tilde M_n(t)=0$. This is carried out by considering first the case $q=2$. 
 By Remark \ref{r1}, it suffices to have $\lim_{n\to\infty} M_n(t)=0$ for $q=2$, and this fact is proved in \cite{dss09}, Proposition 4.3.
 
 To extend the convergence to any $q\in]2,\infty[$,
 we must establish that for some fixed $n_0>0$,
 \beq
 \label{1.50}
 \sup_{n\ge n_0} E\left(\left\Vert\nu_{G_n,Z}^t\right\Vert^q_{L^q(K_a^D(t))}\right) < \infty.
 \eeq
 a result which holds true for $v_{G_n,Z}^t$. Once more, the first step in the proof of \eqref{1.50} consists in obtaining the upper bound
 \beq
 \label{1.51}
 E\left(\left\Vert\nu_{G_n,Z}^t\right\Vert^q_{L^q(K_a^D(t))}\right)\le C E\int_0^t ds \int_{K_a^D(t))} dx E\left(\Vert G_n(t-s,x-\cdot) Z(s,\cdot)\Vert^q_\hac\right).
 \eeq
 This follows easily by applying first Cauchy-Schwarz' inequality to the inner product on $\hact$ and then H\"older's inequality. Once we have \eqref{1.51}, we can obtain \eqref{1.50} by following the steps of
 the proof of Proposition 3.4 in \cite{dss09}.
 
 The last ingredient for the proof of \eqref{t2.2} consist of the extension of Theorem 4.6 in \cite{dss09}. This requires the following additional arguments. Firstly, using similar notations as in that reference, we set
\beqn
R_n^{m,\gamma,D}(t) = E\left(\Vert u_n^{v,(m)}(t)\Vert^q_{W^{\gamma,q}(K_a^D(t))}\right),
\eeqn
where $u_n^{v,(m)}(t,x)$ stands for the $m$-th Picard iteration of a similar equation as \eqref{0.9} with $G$ replaced by the smoothed version $G_n$. In comparison with \cite{dss09}, in order to check that
$\sup_{n,m\ge 1}R_n^{m,\gamma,D} <\infty$, we have to study the additional term
\beqn
T_n^{m,\gamma,D,3}(t) := E\left(\Vert \nu^t_{G_n,\sigma(u_n^{v,(m)})1_{K_a^D})}\Vert^q_{W^{\gamma,q}(K_a^D(t))}\right)
\eeqn 
and more specifically, to check that
\beq
\label{1.52}
T_n^{m,\gamma,D,3}(t)\le C_1+C_2\int_0^t ds \ R_n^{m-1,\gamma,D}(s),
\eeq
for some positive constants $C_1$, $C_2$.

This property holds true when $T_n^{m,\gamma,D,3}(t)$ is replaced by
\beqn 
E\left(\Vert v^t_{G_n,\sigma(u_n^{v,(m)})1_{K_a^D})}\Vert^q_{W^{\gamma,q}(K_a^D(t))}\right)
\eeqn
(see the arguments on page 42 of \cite{dss09} based upon Proposition 3.5 of this reference).
In a similar way, \eqref{1.52} follows from Proposition \ref{p4} and more precisely, from \eqref{1.9}. 

This completes the proof of \eqref{t2.2}.
\medskip

An important consequence of \eqref{t2.2} is the following. For any $t>0$, a.s., the sample paths of
$(u^{\ep,v}(t.x)1_{K^D(t)}(x), x\in\IR^3)$ are $\alpha$-H\"older continuous with $\alpha\in\mathcal{I}$. 
In addition, for any $q\in[2,\infty[$,
\beq
\label{1.53}
\sup_{\ep\in]0,1], v\in\mathcal{P}_\hac^N}\sup_{t\in[0,T]}E\left(\vert u^{\ep,v}(t,x)-u^{\ep,v}(t,y)\vert^q\right) \le C\vert x-y\vert^{\alpha q},
\eeq
for any $x,y\in K^D(t)$, $\alpha\in\mathcal{I}$. Hence, in order to prove \eqref{t2.3} it remains to establish that, for any $q\in[2,\infty[$ and
$\alpha\in\mathcal{I}$, there exists $C>0$ such that for every $t, \bar t\in[0,T]$,
\beq
\label{1.54}
\sup_{\ep\in]0,1], v\in\mathcal{P}_\hac^N}\sup_{x\in D}E\left(\vert u^{\ep,v}(t,x)-u^{\ep,v}(\bar t,x)\vert^q\right) \le C\vert t-\bar t\vert^{\alpha q},
\eeq
For this, we will follow the steps of Section 4.3 in \cite{dss09} devoted to the analysis of the time regularity of the solution to \eqref{0.2} and get an extension of
Theorem 4.10. 

As in the first part of the proof, we consider the case $\ep=1$. The additional required ingredient consists of showing that
\begin{align}
\label{1.55}
&E\left(\left\vert\int_0^t  \langle G(t-s,x-\cdot)\sigma(u^v(s,\cdot))1_{K^D(s)}(\cdot),v\rangle_\hac\ ds\right.\right.\nonumber\\
&\left.\left.- \int_0^{\bar t} \langle G(\bar t-s,x-\cdot)\sigma(u^v(s,\cdot))1_{K^D(s)}(\cdot),v\rangle_\hac\ ds\right\vert^q\right)\nonumber\\
& \le C \vert t-\bar t\vert^{\alpha q},
\end{align}
uniformly in $x\in D$.

Remark that the stochastic process 
\beqn
\{Z(s,y):= \sigma(u^v(s,y))1_{K^D(s)}(y), \ (s,y)\in[0,T]\times \IR^3\},
\eeqn 
satisfies the assumptions of Proposition \ref{p5} with $\mathcal{O}= D$ and arbitrarily large $q$, This fact is proved in Theorem 4.10 in \cite{dss09}. Thus, \eqref{1.55}
follows from that Proposition.

Going through the arguments, it is easy to realize that for $u^{v,\ep}$, we can get uniform estimates in $\ep\in]0,1]$ and $v\in\mathcal{P}_\hac^N$, and therefore \eqref{1.54}
holds true.
This ends the proof of \eqref{t2.3} and of the Theorem.

\hfill\qed
\bigskip
\begin{rem}
\label{r3}
In connection with conclusion (4.8) of Theorem 4.1 in \cite{dss09}, we notice that property (\ref{t2.1}) implies
\beqn
\sup_{\ep\in]0,1], v\in\mathcal{P}_\hac^N}\sup_{t\in[0,T]}E\left(\Vert u^{\ep,v}(t)\Vert^q_{L^q(K^D(t))}\right) < \infty.
\eeqn
\end{rem}

The estimates on increments described in (\ref{2}) are a consequence of  \eqref{t2.3}. Indeed, as has been already mentioned, for any $v\in\mathcal{P}_{\hac}^N$, the stochastic process $V^v$ is the solution to the particular equation (\ref{0.9}) obtained by setting $\ep=0$. 
\bigskip

\begin{prop}
\label{p3} Assume {\bf (H)}. Consider a family $(v^\ep, \ep>0)\subset \mathcal{P}_\hac^N$ and $v\in\mathcal{P}_\hac^N$ such that a.s.,
\beqn
\lim_{\ep\to 0}\Vert v^\ep-v\Vert_w = 0.
\eeqn
Then, for any $(t,x)\in [0,T]\times D$ and any $q\in[2,\infty[$,
\beq
\label{1.5}
\lim_{\ep\to 0}E\left(\vert u^{\ep,v^\ep}(t,x) - V^v(t,x)\vert^q\right)=0.
\eeq
\end{prop}

\noindent{\it Proof.} We write
\beqn
u^{\ep,v^\ep}(t,x) - V^v(t,x)=\sum_{i=1}^4 T_i^\ep(t,x),
\eeqn
with
\begin{align*}
T_1^\ep(t,x)&= \int_0^t \left[G(t-s)*\left(b(u^{\ep,v^\ep}(s,\cdot))-b(V^v(s,\cdot))\right)\right](x)\ ds,\\
T_2^\ep(t,x)&= \int_0^t \left\langle G(t-s,x-\cdot)\left[\sigma(u^{\ep,v^\ep}(s,\cdot))-\sigma(V^v(s,\cdot))\right], v^\ep(s,\cdot)\right\rangle_\hac ds,\\
T_3^\ep(t,x)&= \int_0^t \left\langle G(t-s,x-\cdot)\sigma(V^v(s,\cdot)), v^\ep(s,\cdot)-v(s,\cdot)\right\rangle_\hac ds,\\
T_4^\ep(t,x)&= \sqrt \ep\sum_{k\ge 1}\int_0^t\langle G(t-s,x-\cdot) \sigma(u^{\ep,v^\ep}(s,\cdot)),e_k\rangle_\hac d B_k(s).
\end{align*}
Fix $q\in[2,\infty[$. H\"older's inequality with respect to the measure on $[0,t]\times \IR^3$ given by $G(t-s,dy)ds $, along with the Lipschitz continuity of $b$ yield
\begin{align*}
&E\left(\vert T_1^\ep(t,x)\vert^q\right)\le \left(\int_0^t ds \int_{\IR^3} G(s,dy)\right)^{q-1}\\
&\times \int_0^t  \sup_{(r,z)\in[0,s]\times\IR^3}E\left(\vert u^{\ep,v^\ep}(r,z)-V^v(r,z)\vert^q\right)\left(\int_{\IR^3} G(s,dy)\right)\ ds\\
& \le C  \int_0^t  \sup_{(r,z)\in[0,s]\times\IR^3}E\left(\vert u^{\ep,v^\ep}(r,z)-V^v(r,z)\vert^q\right)\ ds
\end{align*}
To study $T_2^\ep(t,x)$, we apply Cauchy-Schwarz' inequality to the inner product on $\hac$ and then H\"older's inequality with respect to the measure on  $[0,t]\times \IR^3$ given by $\vert\tf G(s)(\xi)\vert^2 ds\ \mu(d\xi)$. Notice that this measure can also we written as $[G(s)\ast\tilde G(s)](x) \varGamma(dx) ds$.
The Lipschitz continuity of $\sigma$ along with (\ref{noise}) and the property
$\sup_{\ep}\Vert v^\ep\Vert_{\hact}\le N$, imply
\begin{align*}
&E\left(\vert T_2^\ep(t,x)\vert^q\right)\le E\left(\int_0^t \left\Vert G(t-s,x-\cdot)\left[\sigma(u^{\ep,v^\ep}(s,\cdot))-\sigma(V^v(s,\cdot))\right]\right\Vert^2_{\hac} ds\right)^{\frac{q}{2}}\\
& \times \left(\int_0^t \Vert v^\ep(s,\cdot)\Vert_{\hac}^2 ds\right)^{\frac{q}{2}}\\
& \le C E\left(\int_0^t \left\Vert G(t-s,x-\cdot)\left[\sigma(u^{\ep,v^\ep}(s,\cdot))-\sigma(V^v(s,\cdot))\right]\right\Vert^2_{\hac} ds\right)^{\frac{q}{2}}\\
& \le C \left(\int_0^t ds \int_{\IR^3} \mu(d\xi)\vert\tf G(t-s)(\xi)\vert^2\right)^{\frac{q}{2}-1}\\
& \times \int_0^t  \sup_{(r,z)\in[0,s]\times\IR^3}E\left(\vert u^{\ep,v^\ep}(r,z)-V^v(r,z)\vert^q\right)\left(\int_{\IR^3} \mu(d\xi)\vert\tf G(s)(\xi)\vert^2\right)\ ds\\
& \le C \int_0^t  \sup_{(r,z)\in[0,s]\times\IR^3}E\left(\vert u^{\ep,v^\ep}(r,z)-V^v(r,z)\vert^q\right)\ ds.
\end{align*}
For any $(t,x)\in[0,T]\times \IR^3$, the stochastic process 
\beqn
\{G(t-s,x-y)\sigma(V^v(s,y)), (s,y)\in[0,T]\times \IR^3\}
\eeqn
 satisfies the property 
\beq
\label{h}
\sup_{v\in \mathcal{P}_{\hac}^N}\sup_{s\in[0,T]} E\left(\left\Vert G(t-s,x-\cdot)\sigma(V^v(s,\cdot))\right\Vert^q_{\hac}\right)<\infty.
\eeq
Indeed, by applying (\ref{t2.1}) to the particular case $\ep=0$, we get
\beq
\label{1.14}
\sup_{v\in\mathcal{P}_{\hac}^N} \sup_{(t,x)\in[0,T]\times \IR^3} E\left(\vert V^v(t,x)\vert^q\right)<\infty.
\eeq
Then, we apply H\"older's inequality with respect to the measure on $\IR^3$ given by $\vert\tf G(t-s)(\xi)\vert^2 \mu(d\xi)$, along with the linear growth property of $\sigma$, and we obtain
\begin{align*}
& E\left(\left\Vert G(t-s,x-\cdot)\sigma(V^v(s,\cdot))\right\Vert^q_{\hac}\right)\nonumber\\
& \le C \left(\int_{\IR^3} \vert\tf G(t-s)(\xi)\vert^2\mu(d\xi)\right)^{\frac{q}{2}}\nonumber\\
& \times\left(1+ \sup_{(s,y)\in[0,T]\times \IR^3}E\left(\left \vert V^v(s,y)\right\vert^q\right)\right).
\end{align*}
With (\ref{noise}) and (\ref{1.14}), we have (\ref{h}).

From \eqref{h}, it follows that  $\{G(t-s,x-y)\sigma(V^v(s,y), (s,y)\in[0,T]\times \IR^3\}$
takes its values in $\hact$, a.s.
Since $\lim_{\ep\to 0}\Vert v^\ep-v\Vert_w = 0$, a.s.,
\beqn
\lim_{\ep\to 0}\left\vert \int_0^t \left\langle G(t-s,x-\cdot)\sigma(V^v(s,\cdot)), v^\ep(s,\cdot)-v(s,\cdot)\right\rangle_\hac ds\right\vert = 0.
\eeqn
Applying (\ref{h}) and bounded convergence, we see that the above convergence takes place in $L^q(\Omega)$ as well. Thus,
\beqn
\lim_{\ep\to 0}E\left(\vert T_3^\ep(t,x)\vert^q\right)=0.
\eeqn
By the $L^q$ estimates of the stochastic integral and \eqref{t2.1}, we have
\begin{align*}
&E\left(\left\vert\sum_{k\ge 1}\int_0^t\langle G(t-s,x-\cdot) \sigma(u^{\ep,v^\ep}(s,\cdot)),e_k\rangle_\hac d B_k(s)\right\vert^q\right)\\
& = E\left(\int_0^t  \Vert G(t-s,x-\cdot) \sigma(u^{\ep,v^\ep}(s,\cdot))\Vert_{\hact}^2\ ds\right)^{\frac{q}{2}}\\
& \le \left(\int_0^t ds \int_{\IR^3} \mu(d\xi)\vert\tf G(t-s)(\xi)\vert^2\right)^{\frac{q}{2}-1}\\
& \times \int_0^t  \left(1+\sup_{(r,z)\in[0,s]\times\IR^3}E\left(\vert u^{\ep,v^\ep}(r,z)\vert^q\right)\right)\left(\int_{\IR^3} \mu(d\xi)\vert\tf G(s)(\xi)\vert^2\right)\ ds\\
& \le C \int_0^t  \left(1+\sup_{(r,z)\in[0,s]\times\IR^3}E\left(\vert u^{\ep,v^\ep}(r,z)\vert^q\right)\right)\ ds\\
& \le C.
\end{align*}
This yields
\beqn
\lim_{\ep\to 0}E\left(\vert T_4^\ep(t,x)\vert^q\right)=0.
\eeqn
We end the proof of the Proposition by applying the usual version of Gronwall's lemma.

Notice that we have actually proved the stronger statement
\beq
\label{1.6}
\lim_{\ep\to 0}\sup_{(t,x)\in[0,T]\times\IR^3}E\left(\vert u^{\ep,v^\ep}(t,x) - V^v(t,x)\vert^p\right)=0.
\eeq

\hfill\qed 

\noindent{\it Proof of Theorem \ref{t1}}. As has been argued, it suffices to check the validity of \eqref{2} and \eqref{p}. These statements follow from Theorem \ref{t2} and Proposition \ref{p3}, respectively.

\hfill\qed

\end{document}